\documentstyle[12pt]{article}
\begin{document}
\begin{center}

{\bf Nikol'skii-type inequalities }\\

\vspace{3mm}

  {\bf for rearrangement invariant spaces }\\

\vspace{3mm}

{\sc Ostrovsky E., Sirota L.}\\

\normalsize

\vspace{3mm}
{\it Department of Mathematics and Statistics, Bar-Ilan University,
59200, Ramat Gan, Israel.}\\
e-mail: \ galo@list.ru \\

{\it Department of Mathematics and Statistics, Bar-Ilan University,
59200, Ramat Gan, Israel.}\\
e - mail: \ sirota@zahav.net.il \\

\end{center}

\vspace{4mm}

  {\it Abstract.} In this paper we generalize the classical Nikol'skii inequality on the many popular classes pairs of rearrangement invariant
(r.i.) spaces  and construct some examples in order to show the exactness
of our estimations.\\



\vspace{3mm}

{\bf 1. Introduction. Notations. Statement of problem.}\\

 Firs of all we recall the classical Nikol'skii inequality:
there exist a values $ p_-,p_+, q_-, q_+, \ 1 \le p_- < p_+ < q_- < q_+ \le \infty, $
such that for all the values $ p \in(p_-, p_+), \ q \in (q_-, q_+) $
the following inequality holds:

 $$
 |t_n|_q \le C \ \sigma^{1/p - 1/q} \ |t_n|_p, \eqno(1)
 $$
  and automatically $ \ 1 \le p \le q \le \infty. $\par

  Here $ t_n \in A(n), $ where $ A(n) $ is the (sub)space of
 trigonometrical or algebraic polynomial from several variables
 defined on the sets $ T = [-\pi, \ \pi]^l, \ T = [0,1]^l $ with general
degree $ \le n, \ n = 1,2,3,\ldots, $ or entire function of exponential type
$ \le n, \ n \in [1,\infty) $ defined on the set $ R^l, \ \sigma = \sigma(n) $ is some increasing function; as a rule $ \sigma = C_1 \ n^{\gamma}, \ \gamma > 0. $ \par

   Recall that the entire exponential function $ t_n(x), \ x \in R^l $ of finite type order $ \le n $ may be defined as follows:

$$
t_n(x) = \int_{ [-n,n]^l } \exp(i x \xi) \ g(\xi) \ d \xi, \ g(\cdot) \in L_2(R^l).
$$

 Hereafter $ C, C_j $ will denote any non-essential finite positive constants. As usually, for the measurable function
 $ f: T \to R, $ where $ (T, M, \mu) $ is a measurable space with non-trivial
 sigma-finite measure $ \mu, $

 $$
|f|_p(T, \mu) = |f|_p(\mu) = |f|_p = \left[\int_T |f(x)|^p \ \mu(dx) \right]^{1/p}, \ p < \infty,
 $$
  $ L_p(T,\mu) = L_p(\mu) = \{f: \ |f|_p < \infty \}; \ m $
 will denote usually Lebesque measure, and we will write $ m(dx) = dx; $
  $ |f|_{\infty} \stackrel{def}{=} \sup_x |f(x)|. $ \par
  There are many generalizations of inequality (1) on the case if $ T $ is convex
  polytop in the space $ R^l, \ $ ( [2], [35], [36], [39]); "weight" Nikol'skii
   inequalities for polynomials defined on the set
  $ [-1, \ 1] $ with Jakobi weight; for polynomial defined on the set $ R^l $ with exponential weight $ w = w(x) $ of a view $ w(x) = $

  $$
   \exp \left(-|x|^{\theta} \right), \ \theta > 0, |x| = |(x_1,x_2,\ldots,x_l)|
= \sqrt{x^2_1 + x^2_2 + \ldots + x^2_l},
  $$

see [39], [40], [41] etc. \par
    There exists a general definition of the {\bf Nikol'skii class }
$ \{ A(n) \} $ over measurable triple $ (T,M,\mu) $ with some function $ \sigma = \sigma(n) $ see [2], where {\it the inequality (1) is postulated} under some
increasing sequence (function) $ \sigma = \sigma(n) $ such that

$$
\lim_{n \to \infty} \sigma(n) = \infty, \ \sigma(n) \ge 3.
$$

 In this definition instead the subspace of trigonometrical or algebraical polynomials with common degree $ \le n $  used some monotonically decreased sequence (or set) $ A(n) $
of linear {\it subspaces} of a space

$$
L (p_-, q_+) = \cap_{p \in (p_-, q_+)} L_p(\mu); \ m > n \ \Rightarrow
$$

$$
 A(n) \subset A(m); \ \forall n \ A(n) \subset L(p_-, q_+) .
$$
  {\bf Further we will consider only the triples  $ (T,M,\mu) $  with
some (non-trivial) Nikol'skii class  $ A(n). $ } \par

   The Nikol'skii inequality play a very important role in the theory of approximation, theory of function of several variables,
   functional analysis (imbedding theorems for Besov spaces). See, for example, [4], [2], [35] etc. \par

  The inequality (1) may be rewritten as follows. Let $ (X, ||\cdot||X) $ be any
  rearrangement invariant (r.i.) space on the set $ T; $ denote by $ \phi(X, \delta) $
  its fundamental function

  $$
  \phi(X,\delta) = \sup_{A, \mu(A) \le \delta} ||I(A)||X, \ I(A)= I(A,x) =
   1, x \in A,
  $$
  $ I(A) = I(A,x) = 0, \ x \notin A.$ We define also for two function r.i. spaces $ (X, \ ||\cdot||X ) $ and $ (Y, \ ||\cdot||) $
  over our set $ T \ $ and for arbitrary finite positive constants $ K_1,\ K_2 $
  the so-called {\sc Nikol'skii two-space functional, briefly: NF functional }
  between the spaces $ X $ and $ Y $
  $$
  W_n(X,Y, K_1, K_2) \stackrel{def}{=} \sup_{t_n \ne 0, t_n \in A(n)} \frac{||t_n||X}
  {\phi(X,K_1/\sigma) } : \frac{||t_n||Y}{\phi(Y,K_2/\sigma)},
  $$
  $ W_n(X,Y) = W_n(X,Y,1,1).$ Then (1) is equivalent to the following inequality:

  $$
  P \in (p_-, p_+), \ q \in (q_-, q_+) \ \Rightarrow \sup_n W_n(L_q, L_p) < \infty). \eqno(2)
  $$

 \vspace{3mm}

{\bf Definition 1.} \par

\vspace{3mm}

{\sc By definition, the {\it pair} of r.i. spaces $ (X, \ ||\cdot||X ) $ and
  $ (Y, \ ||\cdot||) $ over some fixed Nikol'skii class is said to be a {\it (strong) Nikol'skii pair}, write: $ (X,Y) \in Nik, $ if the NF functional between $ X $ and
  $ Y $  is finite:

  $$
  \sup_n W_n(X,Y) < \infty, \eqno(3)
  $$
 and is called a {\it weak Nikol'skii pair}, write $ (X,Y) \in wNik, $
 if for some non-trivial constants } $ K_1, K_2 $

  $$
  \sup_n W_n(X,Y, K_1, K_2) < \infty, \eqno(4)
  $$

    {\bf Our aim is description some pair of r.i. spaces with strong and weak Nikol'skii properties. } \par
    Roughly speaking, we will prove that the most of popular {\it pairs} of r.i. spaces are strong, or at last weak Nikol'skii pairs. \par
 The paper is organized as follows. In the next section we introduce and investigate a {\it new class } of r.i. spaces, so-called moment rearrangement
invariant spaces, briefly, m.r.i. spaces. In the section 3 we formulate and prove the main result of paper for m.r.i. spaces. \par
 In the pilcrow 4 we consider some examples in order to show the exactness of obtained estimations and investigate the low
bounds in our inequalities. In the section 5 we will receive the
 Nikol'skii inequality for (generalized) Zygmynd spaces. \par
 In the last section 6 we investigate the {\it inverse} Nikol'skii inequality
for Lorentz's spaces in order to emphasize  the {\it precision } of
obtained results. \par

\vspace{4mm}

{\bf 2. Auxiliary facts. Moment rearrangement spaces.}\par

\vspace{4mm}

  Let $ (X, \ ||\cdot||X) $ be a r.i. space, where $ X $ is linear subset on the space of all measurable function $ T \to R $ over our measurable space
 $ (T,M,\mu) $ with norm $ ||\cdot||X $ equipped with Nilol'skii class $ A(n). $
\par

\vspace{3mm}

 {\bf Definition 2.} \par

\vspace{3mm}

 {\sc We will say that the space $ X $ with the norm $ ||\cdot||X $ is {\it moment
 rearrangement invariant space,} briefly: m.r.i. space, or
$ X =(X, \ ||\cdot||X) \in m.r.i., $
 if there exist a real constants $ a, b; 1 \le a < b \le \infty, $ and some {\it rearrangement invariant norm } $ < \cdot > $ defined on the space of a real functions defined on the interval $ (a,b), $ non necessary to be finite on all the functions, such that

  $$
  \forall f \in X \ \Rightarrow || f ||X = < \ h(\cdot) \ >, \ h(p) = |f|_p. \eqno(5)
  $$

   We will say that the space $ X $ with the norm $ ||\cdot||X $ is {\it weak moment
 rearrangement space,} briefly, w.m.r.i. space, or $ X =(X, \ ||\cdot||X) \in w.m.r.i.,$
 if there exist a constants $ a, b; 1 \le a < b \le \infty, $ and some {\it functional } $ F, $ defined on the space of a real functions defined on the
interval $ (a,b), $ non necessary to be finite on all the functions, such that }

  $$
  \forall f \in X \ \Rightarrow || f ||X = F( \ h(\cdot) \ ), \ h(p) = |f|_p. \eqno(6)
  $$

    We will write for considered w.m.r.i. and m.r.i. spaces $ (X, \
  ||\cdot||X) $

   $$
       (a, b) \stackrel{def}{=} supp(X),
   $$
   ("moment support"; not necessary to be uniquely defined)
   and define for other such a space $ Y = (Y, \ ||\cdot||Y ) $ with
   $ (c,d) = supp(Y) $
   $$
   supp(X) >> supp(Y),
   $$
    iff $ \min(a,b) > \max(c,d). $ \par
     It is obvious that arbitrary m.r.i. space is r.i. space.\par

   There are many r.i. spaces satisfied the condition (5): exponential Orlicz's  spaces,
  some Martzinkiewitz spaces, interpolation spaces (see [14], [42]).\par

    In the article [13] are introduced the so-called $ G(p,\alpha) $ spaces
  consisted on all the measurable function $ f: T \to R $  with finite norm
   $$
   ||f||_{p,\alpha} = \left[ \int_1^{\infty} \left(\frac{|f|_x}{x^{\alpha}}
    \right)^p \ \nu(dx) \right]^{1/p},
   $$

   where $ \nu $ is some Borelian measure.\par
    Astashkin in [42] proved that the space $ G(p,\alpha) $ in the case
   $ T = [0,1] $ and $ \nu = m $ coincides with the
   Lorentz $ \Lambda_p( \log^{1-p \alpha}(2/s) ) $ space.  Therefore,
     both this spaces are m.r.i. spaces.\par
     Another examples. Recently (see [5], [7], [43], [18], [27] - [29], [31] -
     [34]  etc.)
     appears the so-called Grand Lebesque Spaces $ GLS = G(\psi) =
    G(\psi; a,b) $ spaces consisting on all the measurable functions
     $ f: T \to R $ with finite norms

     $$
     ||f||G(\psi) = \sup_{p \in (a,b)} \left[ |f|_p /\psi(p) \right]. \eqno(7)
     $$

      Here $ \psi(\cdot) $ is some continuous positive on the open interval
    $ (a,b) $ function such that

     $$
     \inf_{p \in (a,b)} \psi(p) > 0, \ \sup_{p \in (a,b)} \psi(p) = \infty.
     $$
      It is evident that $ G(\psi; a,b) $ is m.r.i. space and
      $ supp(G(\psi(a,b)) = (a,b). $\par

  This spaces are used, for example, in the theory of probability [5], [7],
 [8], [9], [10]; theory of PDE [27], [28], functional analysis
 [23], [43], theory of Fourier series, theory of martingales etc.\par

 The spaces $ G(\psi, a,b) $ are non-separable [43], but they satisfy the Fatou property. As long as its Boyd's indices (see for detail definition [1, chapter 1] ) $ \gamma_-, \gamma_+ $ are correspondingly

$$
\gamma_- = 1/b, \ \gamma_+ = 1/a,
$$
(see [43]), we conclude on the basis the main result of article [14] that the spaces
$ G(\psi; a,b) $ are interpolation spaces not only between the spaces
$ [L_1, \ L_{\infty}] $ but between the spaces $ [L_v, \ L_s] $ for all values $ (v,s); v < a, \ s > b.$ \par
  Since for arbitrary real-valued continuous function $ f $ defined on the set $ [0,1] $

  $$
 ||f||C[0,1] = \sup_{t \in [0,1]}|f(t)| = \lim_{p \to \infty } |f|_p = \sup_{p \in [1,\infty) } \ | \ |f|_p \ |,
  $$
  the space $ C[0,1] $ is m.r.i. space with $ supp( C[0,1]) = [1,\infty) $
 or equally, e.g., $ supp (C[0,1]) = [3, \infty). $ \par
   But all the Besov's spaces $ B_{p,s}^r(T) $ are w.m.r.i., but not are
   m.r.i. spaces.\par
    Let us consider now the (generalized) Zygmund's spaces $ L_p \ Log^r L, $ which may be defined as an Orlicz's spaces over some subset of the space
$ R^l $  with non-empty interior  and with $ N - $ Orlicz's function of a
 view

    $$
    \Phi(u) = |u|^p \ \log^r( C + |u|), \ p \ge 1, \ r \ne 0.
    $$

     \vspace{3mm}

     {\bf Lemma 1.} \\
     {\bf 1.} All the spaces $ L_p \ Log^r L $ over real line with measure $ m $
     with condition $ \ r \ne 0 $ are not m.r.i. spaces.\\
     {\bf 2.} If $ r $ is positive and integer, then the spaces $ L_p Log^r L $   are w.m.r.i. space.\\
     {\bf Proof.} {\bf 1.} It is sufficient to consider the case $ T = R^1 $   with
     the measure $ m $ and the case $ p > 1. $ \par
      There exists a function $ f_0 = f_0(x) $ belonging to the space
    $ L_p \ Log^r L: $
  $$
  \int_T |f_0|^p \ \log^r(C + |f_0|) \ dx < \infty,
  $$
  but such that for all sufficiently small values $ \epsilon > 0 $

  $$
  \int_T |f_0|^{p \pm \epsilon} \ dx = \infty
  $$

 in the case $ p > 1 $ and

  $$
  \int_T |f_0|^{p + \epsilon} \ dx = \infty
  $$
  In the case $ p = 1. $ \par

 Therefore, the interval $ (a,b) $ in the definition of m.r.i. spaces does not exists.\par
  The assertion {\bf 2} it follows from the formulae
  $$
  |f|^p \ [\log|f|]^k = d^k |f|^p /dp^k, \ k = 1,2,\ldots.
  $$

  {\bf Lemma 2.} There exists an r.i. space without the w.m.r.i. property.\par
  {\bf Proof.} On the interval $ T = [0,1] $ with usual Lebesque measure $ m $ there exists a function $ f $ with standard normal (Gaussian) distribution.
This implies, for example, that

$$
\int_T \exp(p f(x)) \ dx = \exp \left(0.5 \ p^2 \right), \ p \in R.
$$

There exist a functions $ g: R \to R $ such that the function $ h(x) = g(f(x)) $ which {\it distribution} can not be uniquely
defined by means of all positive moments, for instance, $ h(x) = g(f(x)) = [f(x)]^3 $ or $ g(x) = \exp(f(x)). $ \par
   Let us consider a two such a functions $ f_1 $ and $ f_2 $ with {\it different} distributions, but with at the same moments, for example:

   $$
   \int_T |f_1|^p \ dx = \int_T |f_2|^p \ dx = \int_T [\exp(f)]^p \ dx = \exp(p^2/2), \ p \in R.
   $$

 We choose the (quasi)-concave positive strictly increasing continuous function
$ \theta(\cdot), \ \theta(0+) = 0, $ for which

 $$
 \int_0^{\infty} \theta( m \{x: \ |f_1(x)| > \lambda \} ) \ d\lambda = \infty,
 $$
 but
 $$
 \int_0^{\infty} \theta( m \{ x: \ |f_2(x)| > \lambda \} ) \ d\lambda < \infty.
 $$
 The Lorentz r.i. space $ \Lambda(T, \theta) $ over $ T = [0,1] $ with the function
 $ \theta(\cdot) $ and the classical norm (see [30], chapter 2, section 2)

  $$
  ||f||L(T,\theta) = \int_0^{\infty} \theta( m \{x: \ |f(x)| > \lambda \} ) \ d \lambda
  $$
 is not w.m.r.i. space.\par

\vspace{4mm}

  {\bf 3. Main result. Nikol'skii inequality for the pairs of m.r.i. spaces.}\\

\vspace{3mm}

{\bf Theorem 1.}\par
{\it Let $ (X, ||\cdot||X) $ be any m.r.i. space with support $ supp(X) =
(c,d) $ relatively the auxiliary norm $ < \cdot >, $ and let $ (Y,
||\cdot||Y ) $ be
another m.r.i. space over at the same triple $ (T,M,\mu)$ relatively the
second auxiliary norm $ << \cdot >> $ and with $ supp(Y) = (a,b), $ where
$ (c,d) >> (a,b)\ $ and suppose $ 1 \le p_- = a, \ p_+ = b; \ q_- =c, \
q_+ = d \le \infty.$ \par
  Then the pair of m.r.i. spaces $ (X, ||\cdot||X) $ and $ (Y, ||\cdot||Y) $
is the (strong) Nikol'skii pair: }

$$
\sup_n W(X,Y) = C(X,Y) < \infty. \eqno(8)
$$

{\bf Note } that the restriction
$$
1 \le p_- = a, \ p_+ = b; \ q_- =c, \ q_+
= d \le \infty
$$
is not loss of generality and is necessary for the cases if $ A(n) $ is the
set of trigonometrical or algebraical polynomials.\par

{\bf Proof} is very simple. It follows from the definition of Nikol'skii
Class $ A(n) $ that for arbitrary non-zero element of $ A(n): \ t_n \in A(n) $

$$
\frac{|t_n|_q }{ (1/\sigma)^{1/p} } \le C \frac{ |t_n|_p}{(1/\sigma)^{1/q} }, \eqno(9)
$$
$ q \in (c,d), \ p \in (a,b). $ Taking from the bide-sides of the inequality (9) the norm $ << \cdot >> $ and taking into account the monotonocity of the norm $ << \cdot >> $ and the equality

$$
<< \delta^{1/p} >> = || I(A)||Y = \phi(Y,\delta), \ \mu(A) = \delta,
$$
where $ \delta = 1/\sigma, $ we obtain:

$$
|t_n|_q \ \phi(Y, 1/\sigma) \le C ||t_n||Y /\sigma^{1/q}.
$$

  Calculate from the bide-side the norm $ < \cdot >, $ we find analogously
$$
||t_n||X \ \cdot \phi(Y, 1/\sigma) \le C ||t_n||Y \ \cdot \phi(X, 1/\sigma),
$$
wich is equivalent to (8). \par

\vspace{4mm}

{\bf 4. Examples. Low bounds.} \par

\vspace{4mm}

{\bf 1.} We consider now a very important for applications examples of $ G(\psi) $
spaces. Let $ a = const \ge 1,
b = const \in (a, \infty]; \alpha, \beta = const. $ Assume also that at
$ b < \infty \ \min(\alpha,\beta) \ge 0 $ and denote by $ h $ the (unique)
root of equation
$$
(h-a)^{\alpha} = (b-h)^{\beta}, \ a < h < b;
  \ \zeta(p) = \zeta(a,b; \alpha,\beta; p) =
$$
$$
(p-a)^{\alpha}, \ p \in (a,h);
\ \zeta(a,b; \alpha,\beta;p) = (b-p)^{\beta}, \ p \in [h,b);
$$
and in the case $ b = \infty $ assume that
 $ \alpha \ge 0, \beta < 0; $ denote
by $ h $ the (unique) root of equation
 $ (h-a)^{ \alpha} = h^{ \beta}, h > a; $ define in this case

$$
\zeta(p) = \zeta(a,b;\alpha,\beta;p) = (p-a)^{\alpha}, \ p
\in (a,h); \ p \ge h \ \Rightarrow \zeta(p) = p^{\beta}.
$$

   Here and further
 $ p \in (a,b) \ \Rightarrow \psi(p) \asymp \nu(p) $ denotes that
$$
0 < \inf_{p \in (a,b)} \psi(p)/\nu(p) \le \sup_{p \in (a,b)} \psi(p)/\nu(p)
< \infty.
$$

 The space $ G = G_T = G_T(a,b;\alpha,\beta)=
G(a,b; \alpha,\beta) $ consists by definition on all the measurable functions
$ f: T \to R $ with finite norm:
$$
||f|| G(a,b; \alpha,\beta)= \sup_{p \in (a,b)} |f|_p \cdot \zeta(p).
$$

 On the other words, $ G(a,b; \alpha,\beta)$ is the $ G(\psi; a,b) $ space with
$ \psi(p) = 1/\zeta(p). $ \par
  These spaces was introduced in [5], [43]; and in this article was also calculated its fundamental functions. \par
 We rewrite here only the asymptotical expression for
$ \phi(G(a,b; \alpha,\beta), \ \delta) $
 as $ \delta \to 0+: $

$$
\phi(G(a,b; \alpha_1,\beta_1),\delta) \sim (\beta_1 b^2/e)^{\beta_1} \cdot
\delta^{1/b_1} |\log \delta|^{-\beta_1};
$$

$$
\phi(G(c,\infty; \alpha_2, - \beta_2), \delta) \sim (\beta_2)^{|\beta_2|}
|\log \delta|^{-|\beta_2|},
$$

$ 1 \le a < b < c < \infty; \ \alpha_1, \alpha_2, \beta_1, \beta_2 \in
(0,\infty). $
 We obtain using the theorem 1: if $ t_n \in A(n), $ then

$$
||t_n||G(c,\infty, \alpha_2, -\beta_2) \le C(\alpha_1,\alpha_2,\beta_1,
\beta_2, a, b, c) \times
$$

$$
n^{1/b} \ (\log n)^{- \beta_1 - \beta_2} \ ||t_n||G(a,b, \alpha_1, \beta_1).
$$

\vspace{4mm}

Let now $ X = G(a_1, b_1, \alpha_1, \beta_1), \ Y =
G(a_2, b_2, \alpha_2, \beta_2), $ where
$ 1 \le a_2 < b_2 < a_1 < b_1 < \infty. $ We conclude on the basis of Theorem
1 and the expression for $ \phi(G(\psi), \delta) $ written above: $ t_n \in
A(n), \ n \ge 3 \ \Rightarrow $

$$
||t_n||G(a_1, b_1, \alpha_1, \beta_1) \le C_1 \ n^{1/b_2 - 1/b_1}  \
[\log \ n]^{\beta_1 - \beta_2} \ ||t_n||G(a_2, b_2, \alpha_2, \beta_2),
$$
where
$$
C_1 = C_1(a_1,a_2,b_1,b_2,\alpha_1,\alpha_2, \beta_1, \beta_2).
$$

{\bf 2.} In this subsections we will construct some examples in order to illustrate the exactness of result of section 3.\par
 We consider here only the one-dimensional case $ l = 1 $ and only the "trigonometrical case", i.e. $ T = [-\pi, \pi] $ and $ A(n) $ is the collection of trigonometrical polynomials with degree $ \le n $ ("discrete case")
or $ T = R $ and $ A(n) $ is the set of all entire functions of order
$ \le n; n \in [3, \infty) $ ("continuous case"). \par

{\bf Theorem 2.} {\it Let } $ G(\psi) $ {\it and } $ G(\nu) $ {\it be two arbitrary examples of } $ G(\psi) $ {\it spaces. We assert that}

$$
\underline{\lim}_{n \to \infty} \sup_{t_n \in A(n) } W_n(G(\psi), G(\nu) ) =
C > 0. \eqno(10)
$$

{\bf Note } that we do not suppose here

$$
supp \ G(\psi) >> supp \ G(\nu).
$$

{\bf Proof.}
 Let us consider the following functions:

$$
D_n(x) = \frac{\sin^2(nx/2)}{n^2 \ x^2}, \ x \ne 0; \ D_n(0) = 1/4
$$

in the continuous case, i.e. $ T = R $  and

$$
D_n(x) = \frac{\sin^2(nx/2)}{n^2 \ \sin^2(x/2)},  x \ne 0; \ D_n(0) = 1
$$
in the discrete case $ T = [-\pi, \pi] $ (renormed Fejer's kernels). It is
well-known that $ D_n(\cdot) \in A(n). $ \par
 We restrict itself only by continuous case. Let us calculate all the moments
of $ D_n. $ We have:

$$
\int_{-\infty}^{\infty} |D_n(x)|^p \ dx = \int_{-\infty}^{\infty}
 \left| \ \frac{\sin(n \ x/2)}{n \ x } \right|^{2p} \ dx =
$$

$$
n^{-1} \ 2^{1 \ - \ 2p} \ \int_{-\infty}^{\infty} \left| \frac{\sin y}{y} \right|^{2p}
\ dy = n^{-1} \ 2^{1 \ - \ 2p} \ I(2p),
$$

where

$$
I(s) \stackrel{def}{=} \int_{-\infty}^{\infty} \left| \frac{\sin y}{y} \right|^s
\ dy,
$$

and $ s \in [2,\infty).$ \par

It is evident that at $ s \in [2, 10] \ I(s) \asymp C_1. $ It follows from the
saddle-point (or Laplace) method that as $ s \to \infty \ \Rightarrow $

$$
I(s) \sim \frac{C }{ s^{1/2} }.
$$

 Since

$$
1 = \lim_{s \to \infty} s^{1/s} < \sup_{s \ge 1} s^{1/s} = e^{1/e} < \infty,
$$
we conclude that there exist two finite {\it absolute} positive constants
$ (C_-, C^+), \ C_- < C^+, $ such that for all values $ n \ge 1 $ and
$ p \ge 1 $

$$
C_- \ n^{-1/p} \le \left| D_n \right|_p \le C^+ \ n^{-1/p}.
$$

 We receive calculating the $ G(\psi) $ norm of a function $ D_n(\cdot): $

$$
||D_n||G(\psi) = \sup_{p \ge 1} \left[ |D_n|_p/\psi(p) \right]
\ge C_- \ \sup_p \left[ n^{-1/p} \ /\psi(p) \right] =
$$

$$
C_- \ \phi(G(\psi), 1/n)
$$

and we find analogously

$$
||D_n||G(\psi) \le C^+ \ \phi(G(\psi), 1/n).
$$

We have therefore for two spaces $ G(\psi) $ and $ G(\nu): $

$$
W_n(G(\psi), G(\nu)) \ge C_-/C^+ \stackrel{def}{=} C_3,
$$

Q.E.D. \par

{\bf Corollary. }
 A rouge estimation for the constant $ C_3 $ show that $ C_3 \ge 1/9. $ But
in the source Nikol'skii inequality (1) the exact value $ C $ is equal to 2.\par
 Considering instead $ D_n(\cdot) $ a more general function
$ D_n^{\alpha,\beta}(\cdot) \ \in A(n) $ of a view, for example,

$$
D_n^{\alpha,\beta}(x) =
\int_{-n}^{n} \exp(i t x) \ \left[(1 - (|t|/n)^{\alpha} \right]^{\beta} \ dt,
$$
where $ \alpha, \ \beta = const > 0, $ we get after the optimization over
$ \alpha, \ \beta $ by means of computer computations that in the considered before trigonometrical and one-dimensional case

$$
1.4713\ldots \le \underline{\lim}_{n \to \infty} W_n(G(\psi), G(\nu)) \le
$$

$$
\sup_n W_n(G(\psi), G(\nu)) \le 2.
$$

{\bf 3.} Now we intend to generalize the low estimation on the Orlicz's space
$ L(\Phi), $ where
 $ \Phi(\cdot) $ is some Orlicz's function, i.e. even, twice continuous differentiable,
 convex, $ \lim_{u \to 0} \Phi(u)/u = 0, \ \lim_{u \to \infty} \Phi(u)/u = \infty $ etc. \par

 Recall that the Orlicz's norm of a function $ f: T \to R $ may be calculated by equality:

 $$
 ||f||L(\Phi) = \inf_{v > 0} \left[ \frac{1 + \int_T \Phi(v f(x)) \ dx}{v} \right]
 $$

or up to norm equivalence

$$
||f||L(\Phi) = \inf \left\{v: v > 0, \int_T \Phi(|f(x|/v) dx \le 1 \right\},
$$
see [19], [20], [21], [24]. \par

{\bf Theorem 3.}{\it Let $ \Phi_1(\cdot), \ \Phi_2(\cdot) $ be a two Orlicz's
functions such that for some } $ m = 1,2,3,\ldots $ and all the values
$ v \ge 1 $

$$
\int_1^{\infty} \Phi_i(v/y^{2m}) \ dy \le \Phi_i(C_m v), \ i = 1,2. \eqno(11)
$$

{\it Then there exist a two positive finite constants $ K_1, K_2, K_3 = K_{1,2,3}(\Phi_1, \Phi_2, T) $ such that for Nikol'skii  functional
$ W_n( L(\Phi_1), L(\Phi_2) ) $ there holds: }

$$
\underline{\lim}_{n \to \infty}W_n(L(\Phi_1), L(\Phi_2), K_1, K_2) \ge K_3.
\eqno(12)
$$
\vspace{3mm}
{\bf Note} that the condition (11) is satisfies in the cases if $ \Phi(u) =
|u|^p \ S(|u|), \ N(u) = \exp(|u|^p S(u) ), \ u \ge 2, $
where $ S(u) $ is positive slowly varying as $ u \to \infty $
function (see for definition and properties of slowly varying functions
 [26], chapter 1, sections 1.3 - 1.5.) \par
{\bf Proof.} It is enough to consider only the case $ l = 1 $ and $ T = R.$ Let us choose the introduced function $ D_n(x).$ \par
 In order to estimate the Orlicz's $ ||D_n||L(\Phi), \ \Phi = \Phi_{1,2} $ norm of the function $ D_n, $ we need to estimate the following integral:

 $$
n \cdot J = n \cdot J(n, v, \Phi) \stackrel{def}{=} n \cdot \int_T \Phi( v \ D_n(x) ) \ dx =
 $$

$$
\int_T \Phi \left(v \ \frac{\sin^2 (y/2)}{ y^2 } \right) \ dy.
$$

 Note that at the values $ v: \ |v| \le 1 $ the inequality (11) is satisfied. Further we consider only the values $ v \ge 1, $ as long as the function
$ \Phi $ is even.\par

 The {\it low } estimation for the integral $ J $
 $$
 n \ J(n, v, \Phi) \ge \Phi(C_- \ v)
 $$
is evident; we must to obtain the upper bound at the same manner. We have:

$$
n \ J =
\int_{-\pi}^{\pi} \Phi \left(v \ \frac{\sin^2 (y/2)}{ y^2 } \right) \ dy +
$$

$$
2 \int_{\pi}^{\infty} \Phi \left(v \ \frac{\sin^2 (y/2)}{ y^2 } \right) \ dy =
 I_1 + I_2.
$$

The estimation $ I_1 $ is:

$$
I_1 = \int_{-\pi}^{\pi} \Phi \left(v \ \frac{\sin^2 (y/2)}{ y^2 } \right) \ dy \
\le  2 \pi \Phi (0.25 \ v)
$$

since the function $ v \to \Phi(v) $ is monotonic and $ (\sin^2 y/2)/y^2 \le
1/4. $ \par

 Further,

$$
I_2 \le 2 \int_{\pi}^{\infty} \Phi \left( \frac{v}{y^2} \right) \ dy \le
 \Phi(C \ v)
$$
by virtue of condition (11); the replacing $ 1 \to \pi $ in not essential.\par

Therefore,

$$
C^{-1} \ n \ ||D_n||L(\Phi) \asymp \inf_{v > 0} \left( \frac{1}{v} +
\frac{ \Phi(C_1 \ v) }{n \ v} \right). \eqno(13)
$$

 The asymptotically optimal value $ v_0 $ in the right-side of equality (13)
is attained at

$$
v_0 = C_2 \Phi^{-1}(C_3 \ n),
$$

where $ \Phi^{-1} $ denotes the inverse function to the function $ \Phi $ on the
left-part half-line. \par
 We used  here the known ([19],chapter 2, section 9; [20], chapter 2) the expression
 for the fundamental function of Orlicz's spaces. \par

 Thus,

$$
n \ ||D_n||L(\Phi_i) \asymp \frac{ C(1,i) }{ \Phi_i^{-1}(C_{4,i}/n ) } =
C_{5,i} \ n \ \phi( L(\Phi_i), C_{6,i}/n ).
$$

 Substituting into the expression for $ W_n(X, Y), $ where $ X = L(\Phi_1), \
Y = L(\Phi_2), $ we get to the inequality (12). \par
 The general case, if $ m = 2,3,4, \ldots $ may be provided by the choice

$$
t^{m}_{n}(x) = \left[ D_n(x) \right]^m; \ \left[ D_n(\cdot) \right]^m \in
A(m n).
$$

\vspace{4mm}

{\bf 5. The case of (generalized) Zygmund spaces. Other method.}\par

\vspace{4mm}

 We will suppose in this section that the measure $ \mu $ in the triple
$ (T, M, \mu) $ is atomless and that $ n \ge 3.$ \par

 Recall that the (generalized) Zygmund space

$$
 X = L_q \ (Log \ L)^{\gamma}
$$

over source triple is defined as an Orlicz's space with the Orlicz's function
of a view:

$$
 \Phi(u) = \Phi(q, \gamma; u) = |u|^p \ [\log(C(q,\gamma) + u)]^{\gamma},
$$
where $ C(q, \gamma) $ is sufficiently great constant.\par

Note that the fundamental functions for these spaces are as $ n \to \infty: $

$$
\phi(L_q \ (Log \ L)^{\gamma}, C_1/n) \sim C_2 n^{-1/q} (\log n)^{\gamma/q}.
$$
 Let $ Y $ be another Zygmund's space:

$$
 Y = L_p \ (Log \ L)^{- \beta},
$$
where $ q > p \ge 1 $ (the alternative case is trivial). \par
 Since the function $ u \to |u|^p \ [\log(C(q,\gamma) + u)]^{\gamma} $ satisfies
the $ \Delta_2 \ - $ condition, we can use the assertion of theorem 3:

$$
\underline{\lim}_{n \to \infty} \sup_{t_n \ne 0, t_n \in A(n)} \left[\frac
{||t_n||L_q \ [Log \ L]^{\gamma}}{||t_n||L_p \ [Log \ L]^{- \beta} }\right] \ge
$$

$$
C_L(p,q,\gamma,\beta) \ \sigma^{1/p - 1/q } \ [ \log \ \sigma]^{\gamma/q - \beta/p}.\eqno(14)
$$

 The assertion (14) is the {\it inverse } inequality. We ground now the
{\it direct } Nikol'skii inequality for Zygmund spaces at the same manner
as inequality (14), but {\it only in the cases} $ 1 < p < q, \ \gamma,
\beta \ge 0.$ \par

{\bf Theorem 4.} {\it Let }

$$
 1 < p < q, \ \gamma \ge 0, \beta \ge 0.
$$

{\it We assert that for NF functional of considered spaces the following
inequality is true:}

$$
 \sup_{t_n \ne 0, \ t_n \in A(n)}
  \left[\frac {||t_n||L_q \ [Log \ L]^{\gamma}}{||t_n||L_p \
 [Log \ L]^{-\beta} } \right] \le
$$

$$
C_R(p,q,\gamma,\beta) \ \sigma^{1/p - 1/q} \ [ \log \sigma]^{\gamma/q - \beta/p}. \eqno(15)
$$

{\bf Proof.} Since the cases $ \gamma = 0 $ or $ \beta = 0 $ are simple, we investigate further the possibility $ \gamma > 0, \ \beta > 0. $ \par
 It is proved the article [18] that at least for
$ g \in \cup_n A(n) $ and for arbitrary values $ r > q $

$$
||g||L_q \ [Log \ L]^{\gamma} \le C \left[ \frac{r}{r - q} \right]^{\gamma/r} \
|g|_r \eqno(16a)
$$

and analogously may be proved the 'inverse' inequality: for arbitrary
$ s \in (1,p) $

$$
||g||L_p \ [Log \ L]^{-\beta} \ge C \left[ \frac{s}{p - s} \right]^{- \beta/s} \
|g|_s. \eqno(16b)
$$
 We have for the $ t_n \in A(n) $ combined the definition of Nikol'skii
class and inequalities (16a), (16b):

$$
||t_n||L_q \ [Log \ L]^{\gamma} \cdot \le C \left[ \frac{r}{r - q} \right]^{\gamma/r} \cdot |t_n|_r \le
$$

$$
C \left[ \frac{r}{r - q} \right]^{\gamma/r} \cdot \sigma^{1/s - 1/r} \cdot
|t_n|_s \le
$$

$$
C \left[ \frac{r}{r - q} \right]^{\gamma/r} \ \cdot \ \sigma^{1/s - 1/r} \ \cdot
||t_n||L_p \ [Log \ L]^{-\beta} \  \cdot \left[ \frac{s}{p - s} \right]^{ - \beta/s} =
$$

$$
C \ \cdot \ ||t_n||L_p \ [Log \ L]^{- \beta} \ \times Z,
$$

where

$$
Z = \sigma^{1/s - 1/r} \left[ \frac{r}{r - q} \right]^{\gamma/r} \ \cdot \
\left[ \frac{s}{p - s} \right]^{- \beta/s}, \eqno(17)
$$

$$
1 < s < p < q < r. \eqno(18)
$$

 Minimizing the variable $ Z $ over $ (s,r) $ under the restrictions (18)
for sufficiently greatest  values $ \sigma; \ \sigma \ge 3, $
we prove the desired inequality of theorem 4.\par
 More simple, we can choose in order to prove theorem 4 in the expression
(17) for all sufficiently great values $ n $

$$
r = r_0 = q + \frac{\gamma}{q \ \log \sigma}, \ s = s_0 = p – \
\frac{\beta}{p \ \log \sigma}.
$$

\vspace{4mm}

{\bf 6. Lorentz spaces: inverse Nikol'skii inequalities and regular r.i. spaces.}\par

\vspace{4mm}

 We know that the Lorentz's spaces $ \Lambda \left( Log^{-C_1} (C_2/s) \right) $
are m.r.i. spaces  and, following, they satisfy the conclusion of theorem 1. We intend
to construct in this section some examples of {\it low bounds} in the Nikol'skii
inequalities for Lorentz's spaces.\par
 We consider here as before only the one-dimensional trigonometrical case $ T = R $ and consider some Lorentz's spaces over $ R, \ \Lambda(\phi_i), \ i = 1,2, \ \delta \ge 0, $  where $ \phi, \ \phi_i = \phi_i(\delta)) $ are continuous (quasi-)concave non-negative
 strictly increasing functions,  $ \phi_i(0+) = \phi_i(0) = 0. $ \par
  Denote by $ G(\lambda) $ the following function of distribution:

  $$
  G(\lambda) = m \left[y: \frac{\sin^2(y/2) } {y^2 } > \lambda \right], \
  \lambda \in (0, 1/4).
  $$
   It is obvious that $ G(\cdot) $ is continuous, strictly decreasing,  $ G(1/4 - 0) = 0, \ G(0+) = \infty, \ \delta  \in (0, 1/8) \ \Rightarrow  G(\delta)
   \asymp C \delta^{-1/2}. $ \par

  Let us introduce the following  conditions on the functions
  $ \phi: \ \phi(\cdot) \in Q  $ iff

  $$
  \forall \epsilon > 0 \ \Rightarrow \ \int_0^{1/4} \phi(\epsilon \ G(\lambda)) \ d \lambda  \le \phi(C \ \epsilon) ). \eqno(19)
  $$

   {\bf Note } that the converse inequality to the inequality (19):

$$
  \forall \epsilon > 0 \ \Rightarrow \ \int_0^{1/4} \phi(\epsilon \ G(\lambda))   \ d \lambda \ge \phi(C \ \epsilon) )
  $$

 is always true. \par

\vspace{3mm}

 {\bf Theorem 5.} \\

\vspace{3mm}

 {\bf A.}  If  $ \phi_1 \in Q, \ \phi_2 \in Q $  then
 $ \exists K_1, K_2, K_3 = const \in (0,\infty), $  such that

 $$
 \underline{\lim}_{n \to \infty} W_n(\Lambda(\phi_1), \Lambda(\phi_2), K_1, K_2)   = K_3 > 0.  \eqno(20)
 $$

 {\bf B.}  If  $ \phi \in Q $  and r.i. space $ X $  over $ R $
 is arbitrary
 $ G(\psi) $  space or is the Orlicz's space  $ L(\Phi), $  where the Orlicz's function
 $ \Phi $  satisfies the condition (11), then  $ \exists K_1, K_2 = const \in (0, \infty) \ \Rightarrow $

$$
 \underline{\lim}_{n \to \infty} W_n(\Lambda(\phi), X, K_1, K_2) = K_3 > 0.
 \eqno(21)
 $$

{\bf Proof.} {\bf 1.} We consider as before the function  $ D_n(x).$ Let us estimate
the distribution function for $ D_n. $ We get for the values $ \lambda \in (0, 1/4):$

$$
m \{x: \ D_n(x) > \lambda \} = \int_R I(  D_n(x) > \lambda) \ dx =
$$

$$
n^{-1} \int_R I \left( \frac{\sin^2(y/2)}{y^2}  > \lambda \right) \ d y,
$$
where $ I $ denotes the indicator function.\par
 It is easy  to estimate

$$
\int_R I \left( \frac{\sin^2(y/2)}{y^2} > \lambda \right) \ d y =
m \left \{y: y^{-2} \sin^2(y/2) > \lambda \right \} \asymp
$$

$ G(\lambda), \ \lambda \in (0, 1/4). $ Hence

$$
m \{x: \ D_n(x) > \lambda \} \asymp C n^{-1} \ G(\lambda), \ \lambda \in (0, 1/4)
$$
and

$$
m \{x: \ D_n(x) > \lambda \} = 0
$$

in other case.\par
{\bf 2.} We estimate now the Lorentz norm of a function $ D_n. $ We have based on the
definition of the norm in the Lorentz space:

$$
||D_n|| \Lambda(\phi) = \int_0^{\infty} \phi( m \{x: \ D_n(x) > \lambda \}) \ d \lambda
\asymp
$$

$$
\int_0^{\infty} \phi( C \lambda /n) \ d \lambda.
$$

 The last integral is equivalent, by virtue of condition (19) to $ \phi(C/n). $ \par

{\bf 3.} Let now $ \phi_1 \in Q $ and $ \phi_2 \in Q. $ We conclude
repeating the consideration of section 4 for sufficiently great values $ n $ and
taking into account that the function $ \phi(\cdot) $ is the fundamental function for
$ \Lambda(\phi) $ space: $ W_n(\Lambda(\phi_1), \Lambda(\phi_2), K_1, K_2 ) \ge $

$$
\frac{||D_n|| \Lambda(\phi_1)}{\phi_1(K_1/n) } : \frac{||D_n|| \Lambda(\phi_2)}{\phi_2(K_2/n) } \ge
$$

$$
\frac{\int_0^{\infty} \phi_1(G(\lambda) /n) \ d \lambda }{\phi_1(K_1 /n) } :
 \frac{\int_0^{\infty} \phi_2( G( \lambda) /n) \ d \lambda }{\phi_2(K_2 /n) } \ge C            \eqno(22)
$$

  Therefore

 $$
 \underline{\lim}_{n \to \infty} W_n(\Lambda(\phi), X, K_1, K_2)  > 0.
 $$

{\bf 4.} The last assertion of theorem 5 provided analogously. \\

{\bf Concluding remark.} The r.i. space $ (X, \ || \cdot ||X $ is said to be regular
r.i. space, if

$$
\exists C \in (0, \infty) \ \Rightarrow ||D_n||X \asymp \phi(X, C/n). \eqno(23)
$$

 For example, $ G(\psi) $ spaces, Zygmund, Lorentz spaces are regular r.i. spaces. \par

  We can generalize for two regular r.i. spaces $ X, Y $ over $ R $ and $ A(n) = $ collection of exponential type $ \le n $ there exist a pair of non-trivial
  constants $ K_1, K_2 $ such that

 $$
 \underline{\lim}_{n \to \infty} W_n(X,Y, K_1, K_2) > 0.
 $$

\vspace{6mm}


\begin{center}

{\bf REFERENCES }\\

\end{center}

  1. {\sc Bennet C., Sharpley R.} Interpolation of operators. Orlando, Academic Press
  Inc., (1988).\\

  2. {\sc Ditzan Z., Tikhonov S.} (2005). Ul'yanov and Nikolskii-type inequalities. Journal of Approximation Theory. 133(3) (2005), 100-133.\\

  3. {\sc Belinsky E., Dai F., Ditzian Z.} Multivariate approximating averages.
  Journal of Approximation Theory. 125(1) (2003), 85-1105.\\

   4. {\sc Nikol'skii S.M.} Inequalities for entire analytic functions of finite order
   and their applications to the theory of differentiable functions of several
   variables. Trudy Math. Inst. Steklov, 38 (1951), 244 - 278.\\

   5. {\sc Kozachenko Yu. V., Ostrovsky E.I.} (1985). The Banach Spaces of
      random Variables of subgaussian type. {\it Theory of Probab. and Math.
      Stat.} (in Russian). Kiev, KSU, {\bf 32}, 43 - 57.\\

   6. {\sc Ledoux M., Talagrand M.} (1991) Probability in Banach Spaces.
      Springer, Berlin, MR 1102015.\\

  7. {\sc Ostrovsky E.} Bide-side exponential and moment inequalities
     for tail of distribution of Polynomial Martingales. Electronic
     publication, arXiv: math.PR/0406532 v.1 Jun. 2004. \\

  8. {\sc Ostrovsky E.I.} (1999). Exponential estimations for Random Fields
     and its applications (in Russian). Russia, OINPE.\\

  9. {\sc Ostrovsky E.I.} (2002). Exact exponential estimations for random
     field maximum distribution. {\it Theory Probab. Appl.} {\bf 45} v.3,
      281 - 286. \\

   10. {\sc Talagrand M.} (1996). Majorizing measure: The generic chaining.
      {\it Ann. Probab.} {\bf 24} 1049 - 1103. MR1825156 \\

   11. {\sc Talagrand M.} (2005). {\it The Generic Chaining. Upper and
     Lower Bounds of Stochastic Processes.} Springer, Berlin. MR2133757.\\

   12. {\sc Jawerth B., Milman M.} Extrapolation Theory with Applications.
      Mem. Amer. Math. Soc., 440, (1991)\\

   13. {\sc Lukomsky S.F.} About convergence of Walsh series in the spaces nearest to
     $ L_{\infty}. $ Matem. Zametky, 2001, v.20 B.6,p. 882 - 889.(Russian).\\

   14. {\sc Astashkin S.V.} About interpolation spaces of sum spaces, generated by
      Rademacher system. RAEN, issue MMMIU, 1997, v.1 $ N^o $ 1, p. 8-35.\\

   15. {\sc Nevai, Totik V.} Sharp Nikol'skii-type inequalities with exponential weight.
    Annal. Math., 13 (4),(1987), p. 261 - 267.\\

    16. {\sc Grand R., Santicci P.} Nikol'skii-type and maximal inequalities for generalized trigonometrical polynomials. Manuscripta Math., 99(4), 1999, p. 485 - 507.\\

    17. {\sc Timan M.F.} Orthonormal system satisfies an inequality of S.M.Nikol'skii.
    Annal. Math., 4(1), 1978, p. 75 - 82.\\

   18. Capone C., Fiorenza A., Krbec M. On the Extrapolation Blowups in the
   $ L_p $ Scale. \\

   19. {\sc Krasnoselsky M.A., Rutisky Ya.B.} Convex functions and Orlicz's
    Spaces. P. Noordhoff LTD, The Netherland, 1961, Groningen. \\

   20. {\sc Rao M.M, Ren Z.D.} Theory of Orlicz Spaces. - New York, Basel. Marcel
    Decker, 1991. \\

   21. {\sc Davis H.W., Murray F.J., Weber J.K.} Families of $ L_p - $ spaces with
     inductive and projective topologies. Pacific J.Math. - 1970 - v. 34,
    p. 619 - 638.\\

   22. {\sc Steigenwalt M.S. and While A.J.}  Some function spaces related to
    $ L_p. $  Proc. London Math. Soc. - 1971. - 22, p. 137 - 163.\\

   23. {\sc Ostrovsky E., Sirota L.} Fourier Transforms in Exponential Rearrangement
     Invariant Spaces. Electronic Publ., arXiv:Math., FA/040639, v.1, - 20.6.2004.\\

   24. {\sc Rao M.M., Ren Z.D.} Application of Orlicz Spaces. - New York, Basel.
     Marcel Decker, 2002. \\

    25.{\sc Zygmund A.} Trigonometrical Series. Vol. 1. Cambridge University Press
    (1959)\\

   26. {\sc Seneta E.} Regularly Varying Functions. 1985, Moscow edition.\\

   27. {\it Bongioanni B., Forzani L., Harboure E.} Weak type and restricted weak
      type (p,p) operators in Orlicz spaces. (2002/2003). Real Analysis Exchange,
      Vol. 28(2), pp. 381 - 394.\\

   28. {\sc Harboure E., Salinas O. and Viviani B.}  Orlicz Roundedness for
      Certain Classical Operators. Colloquium Matematicum 2002, 91 (No 2), pp. 263 - 282.\\

   29. {\sc Cianchi A.}  Hardy inequalities in Orlicz Spaces. Trans. AMS,
       {\bf 351}, No 6 (1999), 2456 - 2478.\\

   30. {\sc Krein S.G., Petunin Yu., and Semenov E.M.} Interpolation of linear
        operators. AMS, 1982.\\

   31. {\sc A.Fiorenza.} Duality and reflexivity in grand Lebesgue spaces.
       Collectanea Mathematica (electronic version), {\bf 51}, 2, (2000), 131 - 148.\\

   32. {\sc A. Fiorenza and G.E. Karadzhov.} Grand and small Lebesgue spaces and
       their analogs. Consiglio Nationale Delle Ricerche, Instituto per le
      Applicazioni del Calcoto Mauro Picine", Sezione di Napoli, Rapporto tecnico n.
      272/03, (2005).\\

   33. {\sc T.Iwaniec and C. Sbordone.} On the integrability of the Jacobian under
      minimal hypotheses. Arch. Rat.Mech. Anal., 119, (1992), 129 – 143.\\

   34. {\sc T.Iwaniec, P. Koskela and J. Onninen.} Mapping of finite distortion:
   Monotonicity and Continuity.  Invent. Math. 144 (2001), 507 - 531.\\

   35. {\sc Dai F., Ditzian Z., Tikhonov S.} Sharp Jackson inequalities.
    Journal of Appr. Theory, (2007), 04.015 \\

   36. {\sc Grand R., Santucci P.} Nikol'skii-type and maximal inequality for generalized
    trigonometric polynomials. Manuscripta Math.,  99(4), (1999), 485-507.\\

    37. {\sc Levin E., Lubinsky D.} Orthogonal polynomials for exponential weight.
     Canadian Math. Soc., v. 4 (2001)  44 - 51.\\

     38. {\sc DeVore R.A., Lorentz G.G.} Constructive Approximation. Springer, Berlin, (1993)\\

      39. {\sc Mhaskar H.N.} Intruduction to the Theory of Weighted Polynomial Approximation. World Scientific, Singapore, (1996)\\

      40. {\sc Nessel R., Totic V. } Sharp Nikol'skii inequalities with exponential
       weight. Annal. Math., 13(4), (1987), 261 - 267.\\

       41. {\sc Nessel R., Wilmes G.} Nikol'skii inequalities for trigonometrical polynomials and entire functions of exponential type. J. Austr. Math. Soc.
       Ser. A, (1978), 7 - 18.\\

       42. {\sc Astashkin S.V.} Some new Extrapolation Estimates for the Scale of
         $ L_p \ - $  Spaces. Funct. Anal. and Its Appl., v. 37 $ N^o $ 3 (2003),
         73 - 77. \\

        43. {\sc Ostrovsky E., Sirota L.} Moment Banach Spaces: Theory and Applications.
        HIAT Journal of Science and Engeneering, Nolon, Israel, v. 4, Issue 1 -   2, (2007), 233 - 262.\\

\newpage
\begin{center}

{\bf Ostrovsky E.}\\
\vspace{4mm}

 Address: Ostrovsky E., ISRAEL, 76521, Rehovot, \ Shkolnik street.
5/8. Tel. (972)-8- 945-16-13.\\
\vspace{3mm}
e - mail: {\bf Galo@list.ru}\\

\vspace{6mm}

{\bf Sirota L.}\\
\vspace{3mm}
 Address: Sirota L., ISRAEL, 84105, Ramat Gan.

\vspace{4mm}
e - mail: {\bf sirota@zahav@net.il }\\

\end{center}

\end{document}